\documentclass{mcom-l}

\newtheorem{theorem}{Theorem}[section]
\newtheorem{proposition}[theorem]{Proposition}

\newtheorem{corollary}[theorem]{Corollary}

\theoremstyle{definition}

\theoremstyle{remark}

\numberwithin{equation}{section}
\newcommand{\N}{\mathbb{N}}
\newcommand{\Z}{\mathbb{Z}}
\newcommand{\Q}{\mathbb{Q}}

\begin{document}

\title[Zero-sum free sequences]{Zero-sum free sequences with small sum-set}


\author[G. Bhowmik]{Gautami Bhowmik}
\address{Universit\'e de Lille 1\\
Laboratoire Paul Painlev\'e UMR CNRS 8524\\ 
59655 Villeneuve d'Ascq Cedex\\
France}
\curraddr{}
\email{bhowmik@math.univ-lille1.fr}
\thanks{}

\author[I. Halupczok]{Immanuel Halupczok}
\address{D\'epartement de mathématiques et applications\\
Ecole Normale Sup\'erieure\\
45, rue d'Ulm\\
75230 Paris Cedex 05 -- France}
\email{math@karimmi.de}
\thanks{Supported by the Fondation sciences math\'ematiques de Paris}

\author[J.-C. Schlage-Puchta]{Jan-Christoph Schlage-Puchta}
\address{Mathematisches Institut\\Eckerstr. 1\\79104 Freiburg\\Germany}
\email{jcp@math.uni-freiburg.de}


\date{}

\dedicatory{}

\begin{abstract}
Let $A$ be a zero-sum free subset of $\Z_n$ with $|A|=k$. We compute for
$k\leq 7$ the least possible size of the set of all subset-sums of $A$.
\end{abstract}

\maketitle

\section{Introduction and Results}

For an abelian group $G$ and a subset $B$ of $G$, write $\Sigma(B)
:= \{\sum_{b\in C}b \mid C \subset B\}$ 
for the set of all subset sums of $B$. We say that $B$ is zero-sum free
if it contains no non-empty subset adding up to zero. In this note we
are only interested in finite cyclic groups, and we write
$\Z_n=\Z/n\Z$. 

Define $f_n(k) = \min|\Sigma(B)|-1$, where $B$ runs over all
zero-sum free subsets of $\Z_n$ of cardinality $k$,
and set $f(k) := \min_n f_n(k)$. If there are no zero-sum free sets of
cardinality $k$ in $\Z_n$, we set $f_n(k)=\infty$. This function
arises naturally when 
considering the structure of zero-sum free sequences in $\Z_n$ with
not too many repetitions. For example, Gao and Geroldinger \cite{GGmult}
showed that a sequence $a_1, \ldots, a_m$ in $\Z_n$ with $m>\delta n$
contains a sub-sequence adding up to 0, provided that no element occurs
in the sequence more than $\epsilon n$ times, where $\epsilon$ is a
constant depending only on $\delta$, its precise value being
determined by $f_n$.

The following proposition summarises our knowledge about $f$.
\begin{proposition}
\begin{enumerate}
\item We have $f(1)=1$, $f(2)=3$, $f_n(3)=5$ for $n$
  even, and $f_n(3)=6$ for $n$ odd.
\item We have $f(k)\geq 2k$ for $k\geq 4$, and
  $f(k)\geq\frac{1}{9}k^2$.
\item If $p$ is prime, then $f_p(k)\geq\binom{k+1}{2}-\delta$, where
  $\delta=\begin{cases} 0, & k\equiv 0\;(2)\\ 1, & k\equiv 1\;(2)
\end{cases}$
\end{enumerate}
\end{proposition}
The first statement is straight forward. The second is due to 
Eggleton and Erd\H os \cite{EE} and Olson \cite[Theorem~3.2]{OlMult},
respectively, and the third is due to Olson \cite{Olson}.
In this note we describe the computation of $f_n(k)$ for $k\leq 7$ and
all $n$.

The computations for $k=7$ took 18 hours of CPU-time, the same
algorithm solved $k=6$ within 2 minutes. Hence, even if we only assume
exponential growth of the running time, which appears somewhat
optimistic, the case $k=8$ would require some serious improvements of
the algorithm.

Our main result is the following
\begin{theorem}
We have $f(4)=8$, $f(5)=13$, $f(6)=19$, and $f(7)=24$.
\end{theorem}
In fact, we computed $f_n(k)$ for $k\leq 7$ and all $n$;
the results of these computations are listed in the following table.

How to read the table:

The last column gives an example of a set $B$ of
$k$ elements which has no zero-sum and which has the number
of non-empty sums specified in the third column. 
The second column specifies the conditions on $n$ for this example to work.
Some of the examples of $B$ are only specified for some fixed $n_0$; it is
clear how to turn this into an example for any multiple of $n_0$.

Thus one gets: if the condition in the second column is satisfied, then
$f_n(k)$ has at most the value given in the table. Using a computer we checked
that there are no other examples making $f_n(k)$ smaller.

The boldface values in the third column are the values of $f(k)$.

\newcommand{\immer}{}
\[
\begin{array}{c|c@{\,\,\,}c|c|l}
k & \multicolumn{2}{c|}{\text{cond. on }n} & f_n(k) & \text{Example}\\
\hline
\hline
2
& n \ge 4 & \immer & \mathbf{3} & \{1, 2\} \subset \Z_{n}\\
\hline
3
& n \ge 6 & 2|n & \mathbf{5} & \{1, \frac{1}{2}n, \frac{1}{2}n + 1\} \subset
\Z_{n}\\
& n \ge 7 & \immer & 6 & \{1, 2, 3\} \subset \Z_{n}\\
\hline
4
&& 9|n &  \mathbf{8} & \{3, 1, 4, 7\} \subset \Z_{9}\\
& n \ge 10 & 2|n & 9 & \{1, 2,\frac{1}{2}n, \frac{1}{2}n + 1 \} \subset \Z_{n}\\
& n \ge 12 & 3|n & 9 & \{1, \frac{1}{3}n, \frac{1}{3}n + 1, \frac{2}{3}n + 1 \}
\subset \Z_{n}\\
& n \ge 11 & \immer & 10 & \{1, 2, 3, 4\} \subset \Z_{n}\\
\hline
5
& n \ge 14 & 2|n & \mathbf{13} & \{1, 2,\frac{1}{2}n, \frac{1}{2}n + 1,
\frac{1}{2}n + 2 \} \subset \Z_{n}\\
&& 15|n & 14 & \{-1, 2, 3, 4, 5\} \subset \Z_{15}\\
& n \ge 16 & \immer & 15 & \{1, 2, 3, 4, 5\} \subset \Z_{n}\\
\hline
6
& n \ge 20 & 2|n & \mathbf{19} & \{1, 2, 3, \frac{1}{2}n, \frac{1}{2}n + 1,
\frac{1}{2}n + 2 \} \subset \Z_{n}\\
&& 21|n & 20 & \{-1, 2, 3, 4, 5, 6\} \subset \Z_{21}\\
& n \ge 22 & \immer & 21 & \{1, 2, 3, 4, 5, 6\} \subset \Z_{n}\\
\hline
7
&& 25|n & \mathbf{24} & \{5, 10, 1, 6, 11, 16, 21\} \subset \Z_{25}\\
& n \ge 26 & 2|n & 25 & \{1, 2, 3, \frac{1}{2}n, \frac{1}{2}n + 1, \frac{1}{2}n + 2,
\frac{1}{2}n + 3 \} \subset \Z_{n}\\
&& 27|n & 26 & \{1, -2, 3, 4, 5, 6, 7\} \subset \Z_{27}\\
&n \ge 30 & 3|n & 27 & \{1, 2, \frac{1}{3}n, \frac{1}{3}n + 1, \frac{1}{3}n + 2,
\frac{2}{3}n + 1, \frac{2}{3}n + 2\} \subset \Z_{n}\\
& n \ge 30 & 5|n & 27 & \{1, \frac{1}{5}n, \frac{1}{5}n + 1, \frac{2}{5}n,
\frac{2}{5}n+1,\frac{3}{5}n+1,\frac{4}{5}n+1\} \subset \Z_{n}\\
& n \ge 29 & \immer & 28 & \{1, 2, 3, 4, 5, 6, 7\} \subset \Z_{n}\\
\end{array}
\]

It is clear that $f_n(k)$ is either $\infty$ or less than $n$, so
in particular $f_n(k) = \infty$ if $n \le f(k)$. On the other hand,
for any $n > f(k)$ the table does give an example which yields
$f_n(k) < \infty$. Thus we get:
\begin{corollary}
If $k \le 7$, then $f_n(k) = \infty$ if and only if $n \le f(k)$.
\end{corollary}

\section{Description of the algorithm}

\newcommand{\sys}{\mathcal{E}}
\newcommand{\syseq}{\sys(\sim)}
\newcommand{\lmax}{\ell_{\mathrm{max}}}

Let $k$ and $\ell$ be fixed. We want to check whether there exists
a number $n$ and a zero-sum free set $B = \{b_1, \dots b_k\}
\subset \Z_n$ consisting of $k$ distinct elements such that
$|\Sigma(B)|-1 = \ell$. First we describe how to turn the problem into
an algorithmically decidable one, and then we shall describe how to reduce
the amount of computation so as to solve the problem in real time.

Suppose there exist such $n$ and $B$. Then we get an
equivalence relation $\sim$ on the set of non-empty subsets of
$\{1, \dots, k\}$, defined by
$C \sim C' \iff \sum_{i\in C} b_i =\sum_{i\in C'} b_i$, and this equivalence
relation has $\ell$ equivalence classes. Moreover, the elements $b_i$ of $B$
form a solution modulo $n$ of the system of equations and
inequations $\syseq$, which we define as follows:

\begin{enumerate}
\item
For each $i \ne j$, take the inequation $x_i \ne x_j$.
\item
For each $C \subset \{1, \dots, k\}$, $C \ne \emptyset$,
take the inequation $\sum_{i\in C} x_i \ne 0$.
\item
For each pair $C, C' \subset \{1, \dots, k\}$, take
$\sum_{i\in C} x_i =\sum_{i\in C'} x_i$ or 
$\sum_{i\in C} x_i \ne\sum_{i\in C'} x_i$,
depending on whether $C \sim C'$ or not.
\end{enumerate}

On the other hand, any solution modulo $n$ of this system $\syseq$
defines a set $B$ solving the original problem.

The algorithm now does the following. It iterates through
all possible equivalence relations with at most $\lmax$
equivalence classes, where $\lmax = \frac{k(k+1)}{2} - 1$ is one less than the trivial
upper bound for $f_n(k)$. For each relation $\sim$,
it checks whether there is an inequation $L \ne R$ in $\syseq$
such that $L = R$ lies in the $\Z$-lattice generated by
the equations in $\syseq$. If this is the case, then $\syseq$ is not
solvable modulo $n$ for any $n$.
Otherwise, call $\sim$ an almost-example.

Note that not all almost-examples really yield an example of a set $B$.
For example,
$\sys = \{x_1 \ne 0, x_2 \ne 0, x_1 \ne x_2, 2x_1 = 2x_2 = 0\}$
has no solution modulo any $n$, but none of $x_1 = 0$, $x_2 = 0$, $x_1 = x_2$
lies in $\langle 2x_1 = 0, 2x_2 = 0\rangle_\Z$. In fact, by leaving out any
of the inequations the system becomes solvable for any even $n$. 

For each almost-example, the algorithm solves the system of equations
in $\syseq$ in $\Q/\Z$, that is, we solve it in $\Q$ in the usual way,
but whenever an equation $L=0$ is divided by an integer $a$,
one has to separately consider the cases $\frac{L}{a} = 0$,
$\frac{L}{a} = 1$, \dots, $\frac{L}{a} = a - 1$.
Thus one gets a whole list of solutions, and each solution consists
of a list of variables $x_i$ which can be chosen freely, and
linear expressions for the remaining ones.

Now the algorithm symbolically plugs these linear expressions into the
inequalities $L\ne R$ of $\syseq$. If one gets identically $L=R$, then
this is not a solution of $\syseq$; if however one does
not get identically $L=R$ for any of the inequalities,
then almost all values in $\Q/\Z$ for the free
variables $x_i$ yield a solution of $\syseq$ in $\Q/\Z$, which means
that by multiplying by appropriate $n$, we find solutions
in $\Z_n$.
Finally, the computer prints all those solutions of $\syseq$ in $\Q/\Z$, and
we manually check the necessary conditions on $n$ to make
the example work.

\medskip

The problem is now finite, but the number of equivalence
relations which we have to try is of magnitude the number of equivalence
relations on a set of $2^k$ elements, that is, even for $k=4$ we would
have to check about $10^{10}$ cases. Since each single case requires a
considerable amount of computation, this would already stretch our resources.

We first remark that it turned out that the number of almost-examples
is very small. For example, for $k=7$ and $\ell \le 27$ there are only
$19$ up to permutation of the set $B$, so there is no need to optimise any
part of the algorithm treating the almost-examples. And even though the
search for almost-examples finds some almost-examples in several different
shapes given by permutations of the basic set, removing duplicates takes
almost no time compared to the main search. Thus in the sequel,
we describe how the algorithm works in practice, but we only deal with 
the part searching for almost-examples.

The program starts with a system $\sys$ consisting
only of the inequations (1) and (2). Then it recursively adds
equations and inequations of the form (3) to $\sys$.
As soon as $\sys$ gets inconsistent, the program stops in this branch.
By ``inconsistent'' we mean, as described above, that there exists an
inequation whose negation lies in the lattice generated by the equations.

The program also stops if it can easily prove that the final
equivalence relation will have more than $\lmax$ equivalence classes.
To this end, it searches for a maximal anti-clique using a greedy approach.
Start with an empty anti-clique $\mathcal{A}$.
Iterate through all subsets $C \subset \{1,\dots,k\}$. If for all
$C' \in \mathcal{A}$, $\sys_c$ is inconsistent with the equation
$\sum_{i\in C} x_i =\sum_{i\in C'} x_i$, then add $C$ to $\mathcal{A}$.
The cardinality of $\mathcal{A}$ is a lower bound for the number
of different sums we will finally get.
Whether this method yields good bounds heavily depends on the
order in which the subsets $C$ are considered. We will describe the
order below.

We can greatly reduce the computation time by exploiting symmetry
coming from permutations of the elements of $B$. We use the following
general method: we choose a totally ordered set $\Gamma$, and
for each complete system $\sys$
a function $v_{\sys}\colon \{1,\dots, k\} \to \Gamma$
such that $v_{\sigma(\sys)}(\sigma(i)) = v_{\sys}(i)$
for any permutation $\sigma \in S_k$. We may then restrict
our search to those $\sys$ for which $v_{\sys}$ is (weakly) increasing.
During the computation, the program computes lower and upper bounds for
the values of $v_{\sys}$ and stops if $v_{\sys}$ can not be increasing anymore.

The function $v_{\sys}$ which we use is:
\[
\begin{aligned}
v_{\sys}(i) = (&\text{number of equations in $\sys$ of the form $x_i=a+b$},\\
&
\text{number of equations in $\sys$ of the form $a=x_i+b$},\\
&
\text{number of equations in $\sys$ of the form $x_i+a=b+c$}) \in \N^3
.
\end{aligned}
\]
We use the lexicographical order on these tuples, where the first
entry is the most significant one.

It is important to choose a good order in which to try to add
equations and inequations during the recursion, so that we get
contradictions as early as possible.
A good approach is to start with equations between one and
two element sums: on the one hand, such equations
imply a lot of other equations. On the other
hand, if only few such equations exist, then we already get
a lot of different sums,
which is also helpful. Therefore, the program first treats
the equations of the form $x_i = x_{i'} + x_{i''}$,
then the equations of the form $x_i + x_{i'} = x_{i''} + x_{i'''}$,
and the remaining ones only afterwards. Another advantage of this order is
that we are able to apply the symmetry conditions early.

Now we can explain the order in which the above anti-clique $\mathcal{A}$
is built: as we expect to have a lot of inequations between one
and two element sums, the program tries these sums first when
constructing the anti-clique.

Finally, we mention some of the data structures used to work more
efficiently with $\sys_c$.
\begin{itemize}
\item
Do not perform any consistency check of $\sys_c$ with another equation
or inequation twice: always store the old results.
\item
Keep track of the equivalence relation defined by the
equations which we already added to $\sys_c$.
If $C$ and $C'$ are equivalent, then concerning consistency
of $\sys_c$ we do not need to distinguish between $C$ and $C'$,
i.e.\ we are able to use remembered consistency results more often.

Also update the equivalence relation when we accidentally stumble over
an equation which follows from $\sys_c$.
\item
Each time an equation is added to $\sys$, immediately
put the system of equations into upper triangular form.
\end{itemize}
\bibliographystyle{amsplain}

\end{document}